\documentclass[12pt, reqno]{amsart}

\usepackage{amsmath, amsthm, amscd, amsfonts, amssymb, graphicx, color}
\usepackage[bookmarksnumbered, colorlinks, plainpages]{hyperref}% this package creats hyperlinks.
\usepackage{mathtools}
\usepackage{enumerate}% for creating numbered list.
\usepackage{amscd} % This code allows you to draw the commutative diagrams
\usepackage{xy}
\xyoption{all}
\usepackage{float}

\newtheorem{theorem}{Theorem}[section]

\newtheorem{proposition}[theorem]{Proposition}

\theoremstyle{definition}
\newtheorem{definition}[theorem]{Definition}
\newtheorem{example}[theorem]{Example}

\newtheorem{remark}[theorem]{Remark}
\begin{document}
%-----------------------------------------------PAPER / ARTICLE FORMATINGS---------------------------------------------------
\setcounter{page}{1}

\title[ ]{On Gauss Calendar Formula: An Application of Modular Arithmetic}

\author[K. Nantomah]{Kwara Nantomah}

\address{Department of Mathematics, School of Mathematical Sciences, C. K. Tedam University of Technology and Applied Sciences, P. O. Box 24, Navrongo, Upper-East Region, Ghana. }
\email{\textcolor[rgb]{0.00,0.00,0.84}{ knantomah@cktutas.edu.gh}}

%\dedicatory{This paper is dedicated to Professor ABCD}



%\date{Received: xxxxxx; Revised: xxxxxx; Accepted: xxxxxx.
%\newline \indent $^{*}$ Corresponding author}
% }

\begin{abstract}
In this recreative piece of work, we present Gauss' calendar formula with some examples to demonstrate how it is applied. Then, based on it, we give a formula for determining dates of particular week days of a given month, and some examples are also given to demonstrate how it is applied. The key mathematical concept behind the formulas is modular arithmetic.
\end{abstract} \maketitle

%--------------------------------------------------------BODY OF ARTICLE --------------------------------------------------------

%%%%%%%%%%%%%%%%%%%%%%%%%%%%%%%%%%%%%%%%%%%%%%%%%%%%
\section{Introduction}
%%%%%%%%%%%%%%%%%%%%%%%%%%%%%%%%%%%%%%%%%%%%%%%%%%%%

The calendar system that is largely used in the world today is called the Gregorian calendar and it was established in 1582 by Pope Gregory XIII. Based on this calendar, Carl Friedrich Gauss established a formula for finding the day of the week for any given date. The formula is an application of modular arithmetic which is taught at basic level mathematics. It is stated in \cite[p. 110]{Kraitchik-1953-MR} without a proof. In \cite{Schwerdtfeger-2001}, Schwerdtfeger provided an explanation and proof of the formula as well as its mathematical background. For further information on the subject, the reader may consult the works \cite{Schwerdtfeger-2001} and \cite{Schwerdtfeger-2023}.

In this note, we present Gauss' formula with some examples. Then, based on it, we give a formula  for determining dates of particular week days of a month with some given examples. These formulas will come in handy in situations where the computer calendar is incapacitated. They also serve as powerful tools for mental calculations of day of the week for any given date. The outcomes are restricted to the Gregorian calendar.

%%%%%%%%%%%%%%%%%%%%%%%%%%%%%%%%%%%%%%%%%%%%%%%%%%%%
\section{A Brief History About the Calendar}
%%%%%%%%%%%%%%%%%%%%%%%%%%%%%%%%%%%%%%%%%%%%%%%%%%%%

The story of the calendar dates back to Romulus, the founder and first king of Rome. He introduced a year system of 300 days consisting of ten months. Later, his successor, Numa Pompilius added two months and this was used until Julius Caesar introduced a more accurate year of 365.25 days in 46 BC. The extra 0.25 days was catered for by letting the length of \textit{ordinary year} to be 365 days whilst the length of every fourth year, \textit{a leap year} is 366 days. This is referred to as the Julian calendar. The Julian calendar was widely used until 1582 when there was an accumulated error of 10 days. Then Pope Gregory XIII instituted a second calendar reform to compensate for the error. This was achieved by requiring that, a century year should not be a leap year unless the century number is divisible by 4. In other words, years divisible by 100 are not leap years unless they are divisible by 400 as well. For example, 1700, 1800, 1900, 2100, 2200 and, 2300 are not leap years whilst  1600, 2000, 2400, 2800 and 3200 are leap years. Thus, the calendar we use today is called the Gregorian calendar and it took effect from 15th October 1582.

%%%%%%%%%%%%%%%%%%%%%%%%%%%%%%%%%%%%%%%%%%%%%%%%%%%%
\section{Gauss Calendar Formula}
%%%%%%%%%%%%%%%%%%%%%%%%%%%%%%%%%%%%%%%%%%%%%%%%%%%%

In this section, we present Gauss' calendar formula with some given examples to demonstrate how it is applied.

%------------------- Proposition -------------------------------------
\begin{definition}[Leap Year]
A leap year is a year that is divisible by 4 except for century years that must be divisible by 400 as well. By this definition, 1704, 1824, 1908, 2000, 2024, 3600 are leap years whilst  1500, 1998, 2015, 2025, 2500 are not.
\end{definition}

%------------------- Proposition -------------------------------------
\begin{proposition}\label{prop:Gauss-DF}
Gauss formula for calculating day of the week ($W$) is given by \cite[p. 110]{Kraitchik-1953-MR}
\begin{equation}\label{eqn:GDF-1}
W=(D+M+C+Y)\mod 7
\end{equation}
where $D$ is the date, $M$ is the month, $C$ is the century and $Y$ is the year. If the year is a leap and the month is January or February, then the formula is given by
\begin{equation}\label{eqn:GDF-2}
W=[(D+M+C+Y)-1]\mod 7 .
\end{equation}
\end{proposition}

The four digit year is divided into two parts. The first two digits is denoted by $C$ and the last two digits is denoted by $Y$. The codes for $W$, $M$, $C$ and $Y$ are repectively given in Tables \ref{tab:Code-W}, \ref{tab:Code-M}, \ref{tab:Code-C} and \ref{tab:Code-Y} as defined in \cite[p. 110-111]{Kraitchik-1953-MR}.

\begin{table}[H]
   \caption{Codes for Days of the Week} 
   \label{tab:Code-W}
\begin{center}
\begin{tabular}{ |c|c|c|c|c|c|c|c| }
\hline
\textbf{Code($W$)} & 0 & 1 & 2 & 3 & 4 & 5 & 6 \\ 
\hline
\textbf{Day of Week} & Sat& Sun & Mon & Tue & Wed & Thu & Fri \\ 
\hline
\end{tabular}
\end{center}
\end{table}

\begin{table}[h]
   \caption{Codes for Months} 
   \label{tab:Code-M}
\begin{center}
\begin{tabular}{ | c| c | c |c|c|c|c|c|}
\hline
\textbf{Code($M$)} & 0 & 1 & 2 & 3 & 4 & 5 & 6 \\ 
\hline
\textbf{Month} & January& May & August & February & June& September & April\\ 
          & October&         &             & March       &         & December & July \\ 
          &              &          &            & November &        &                  &        \\ 
\hline
\end{tabular}
\end{center}
\end{table}

\begin{table}[H]
   \caption{Codes for the Century (i.e. First Two Digits of the Year)} 
   \label{tab:Code-C}
\begin{center}
\begin{tabular}{ |c|c| }
\hline
\textbf{Century} & \textbf{Code($C$)} \\ 
\hline
15, 19, 23, 27, \dots  &  1 \\ 
\hline
16, 20, 24, 28, \dots &  0 \\ 
\hline
17, 21, 25, 29, \dots &  5 \\ 
\hline
18, 22, 26, 30, \dots &  3 \\ 
\hline
\end{tabular}
\end{center}
\end{table}

%% -------------- This Table is up 28 Digits---------------------------------------------
%\begin{table}[h]
%   \caption{Codes for Y (i.e. Last Two Digits of the Year)} 
%   \label{tab:Code-Y}
%\begin{center}
%\begin{tabular}{ |c|c|c|c|c|c|}
%\hline
%\textbf{Year} & \textbf{Code ($Y$)} & \textbf{Year} & \textbf{Code ($Y$)} & \textbf{Year} & \textbf{Code ($Y$)} \\ 
%\hline
%00& 0&11&6&22& 6 \\ 
%\hline
%01&1&12&1&23& 0 \\ 
%\hline
%02&2&13&2&24& 2 \\ 
%\hline
%03& 3&14&3&25& 3 \\ 
%\hline
%04&5 &15&4&26& 4 \\ 
%\hline
%05& 6&16&6&27& 5 \\ 
%\hline
%06& 0&17&0&28& 0 \\ 
%\hline
%07&1 &18&1& & \\ 
%\hline
%08& 3&19&2& &  \\ 
%\hline
%09& 4&20&4& &  \\ 
%\hline
%10&5 &21&5& &  \\ 
%\hline
%\end{tabular}
%\end{center}
%\end{table}

\begin{table}[H]
   \caption{Codes for Y (i.e. Last Two Digits of the Year)} 
   \label{tab:Code-Y}
\begin{center}
\begin{tabular}{ |c|c|c|c|c|c|c|c|}
\hline
\textbf{Year} & \textbf{Code($Y$)} & \textbf{Year} & \textbf{Code($Y$)} & \textbf{Year} & \textbf{Code($Y$)} & \textbf{Year} & \textbf{Code($Y$)}\\ 
\hline
00&0&26&4&52&2&78&6 \\ 
\hline
01&1&27&5&53&3&79&0 \\ 
\hline
02&2&28&0&54&4&80&2 \\ 
\hline
03&3&29&1&55&5&81&3\\ 
\hline
04&5&30&2&56&0&82&4 \\ 
\hline
05&6&31&3&57&1&83&5 \\ 
\hline
06&0&32&5&58&2&84& 0\\ 
\hline
07&1&33&6&59&3&85&1 \\ 
\hline
08&3&34&0&60&5&86&2 \\ 
\hline
09&4&35&1&61&6&87&3 \\ 
\hline
10&5&36&3&62&0&88&5 \\ 
\hline
11&6&37&4&63&1&89&6 \\ 
\hline
12&1&38&5&64&3&90&0 \\ 
\hline
13&2&39&6&65&4&91&1 \\ 
\hline
14&3&40&1&66&5&92&3 \\ 
\hline
15&4&41&2&67&6&93&4 \\ 
\hline
16&6&42&3&68&1&94&5 \\ 
\hline
17&0&43&4&69&2&95&6 \\ 
\hline
18&1&44&6&70&3&96&1 \\ 
\hline
19&2&45&0&71&4&97&2 \\ 
\hline
20&4&46&1&72&6&98&3 \\ 
\hline
21&5&47&2&73&0&99&4 \\ 
\hline
22&6&48&4&74&1&& \\ 
\hline
23&0&49&5&75&2&& \\ 
\hline
24&2&50&6&76&4&& \\ 
\hline
25&3&51&0&77&5&& \\ 
\hline
\end{tabular}
\end{center}
\end{table}

In the following remarks, we attempt to interpret how the codes in Table \ref{tab:Code-W}, Table \ref{tab:Code-M}, Table \ref{tab:Code-C} and Table \ref{tab:Code-Y} are obtained.

%------------------- Remark -------------------------------------
\begin{remark}
The codes in Table \ref{tab:Code-W} are obtained as follows. It is assumed that the week starts on Sunday and ends on Saturday.  So the codes in the table are generated by using the formula 
\begin{equation}
\text{Code}(W)=P \mod 7
\end{equation}
where $P$ is the position of the day of the week. For example, Saturday is the 7th day of the week. So the code for Saturday is $W=(7\mod 7)=0$.
\end{remark}

%------------------- Remark -------------------------------------
\begin{remark}
The codes in Table \ref{tab:Code-M} are generated for an ordinary (non-leap) year by the formula \cite{Schwerdtfeger-2001}
\begin{equation}
\text{Code} (M)=T \mod 7
\end{equation}
where $T$ is the total number of days of previous months, starting with January. Thus, January is $(0\mod7)=0$, February is $(31\mod7)=3$, March is $(59\mod7)=3$, April is $(90\mod7)=6$, and so on. For a leap year, the codes for January and February are reduced by 1. 

Also, it is observed that, for an ordinary year,  the code for each month can be interpreted as the number of \textit{week days} from the week day of 1st of January to the week day of 1st of that month. For the case of a leap year, the week days for 1st January and 1st February are reduced by one day. 
\end{remark}

%------------------- Remark -------------------------------------
\begin{remark}
The codes in Table \ref{tab:Code-C} are generated by the formula 
\begin{equation}
\text{Code} (C)=\left[ 2(3-C\text{mod}4) + 1\right] \mod 7
\end{equation}
where $C$ is the first two digits of the four-digit year. 
\end{remark}

%------------------- Remark -------------------------------------
\begin{remark}
The codes in Table \ref{tab:Code-Y} are generated by the formula 
\begin{equation}
\text{Code} (Y)=\left( Y+ \lfloor Y/4 \rfloor \right) \mod 7
\end{equation}
where $ \lfloor . \rfloor$ is the floor function and $Y$ is the last two digits of the four-digit year. 
\end{remark}

We now provide some examples on how the formula is used.

\begin{example}
Which day of the week was 25th December 1911? \\
\textit{Solution}: \\
$D$ (the Date) is 25, \\
$M$ (the code for the Month: December) is 5,\\
$C$ (the code for the Century: 19) is 1, \\
$Y$ (the code for the Year: 11) is 6. \\
Therefore the day of the week is calculated as
\begin{equation*}
W=(25+5+1+6)\mod 7=2
\end{equation*}
which  is Monday.
\end{example}

\begin{example}
Which day of the week was 13th May 1693? \\
\textit{Solution}: \\
$D$ (the Date) is 13, \\
$M$ (the code for the Month: May) is 1,\\
$C$ (the code for the Century: 16) is 0, \\
$Y$ (the code for the Year: 93) is 4. \\
Therefore the day of the week is 
\begin{equation*}
W=(13+1+0+4)\mod 7=4
\end{equation*}
which  is Wednesday.
\end{example}

\begin{example}
Which day of the week was 18th January 1743? \\
\textit{Solution}: \\
$D$ (the Date) is 18, \\
$M$ (the code for the Month: January) is 0,\\
$C$ (the code for the Century: 17) is 5, \\
$Y$ (the code for the Year: 43) is 4. \\
Therefore the day of the week is 
\begin{equation*}
W=(18+0+5+4)\mod 7=6
\end{equation*}
which  is Friday.
\end{example}

\begin{example}
Which day of the week was 23rd February 2004? \\
\textit{Solution}: \\
$D$ (the Date) is 23, \\
$M$ (the code for the Month: February) is 3,\\
$C$ (the code for the Century: 20) is 0, \\
$Y$ (the code for the Year: 04) is 5. \\
Since 2004 is a leap year and the month is February, the day of the week is 
\begin{equation*}
W=[(23+3+0+5)-1]\mod 7=2
\end{equation*}
which  is Monday.
\end{example}

\begin{example}
Which day of the week will be 29th January 2048? \\
\textit{Solution}: \\
$D$ (the Date) is 29, \\
$M$ (the code for the Month: January) is 0,\\
$C$ (the code for the Century: 20) is 0, \\
$Y$ (the code for the Year: 48) is 4. \\
Since 2048 is a leap year and the month is January, the day of the week is 
\begin{equation*}
W=[(29+0+0+4)-1]\mod 7=4
\end{equation*}
which  is Wednesday.
\end{example}

%%%%%%%%%%%%%%%%%%%%%%%%%%%%%%%%%%%%%%%%%%%%%%%%%%%%
\section{Formula for Dates of Particular Week Days of a Month }
%%%%%%%%%%%%%%%%%%%%%%%%%%%%%%%%%%%%%%%%%%%%%%%%%%%%

In this section, we present a formula for determining dates of particular week days of a month with some given examples to demonstrate how it is applied.

%------------------- Proposition -------------------------------------
\begin{proposition}\label{prop:Nantomah-DF}
The formula for determining dates of particular week days of a month is given as
\begin{equation}\label{eqn:Nant-DF-1}
D=[W-(M+C+Y)]\mod 7
\end{equation}
where $W$ is the week day, $M$ is the month, $C$ is the century and $Y$ is the year. If the year is a leap and the month is January or February, then the formula is given by
\begin{equation}\label{eqn:Nant-DF-2}
D=[W-(M+C+Y)+1]\mod 7.
\end{equation}
\end{proposition}

%------------------- Remark -------------------------------------
\begin{remark}
The $D$ in Proposition \ref{prop:Nantomah-DF} gives the first date of the particular week day. The other dates are obtained by adding $7$ to $D$ until you get the last date for the month.
\end{remark}

\begin{example}
Find all Saturdays in December 1998. \\
\textit{Solution}: \\
$W$ (the code for Saturday) is 0,\\
$M$ (the code for the December) is 5,\\
$C$ (the code for the Century: 19) is 1, \\
$Y$ (the code for the Year: 98) is 3. \\
The date of the first Saturday is 
\begin{equation*}
D=[0-(5+1+3)]\mod 7=(-9)\mod 7=5.
\end{equation*}
Therefore the Saturdays in December 1998 are: \textbf{5, 12, 19} and \textbf{26}.
\end{example}

\begin{example}
Find all Fridays in July 1718. \\
\textit{Solution}: \\
$W$ (the code for Friday) is 6,\\
$M$ (the code for the July) is 6,\\
$C$ (the code for the Century: 17) is 5, \\
$Y$ (the code for the Year: 18) is 1. \\
The date of the first Friday is 
\begin{equation*}
D=[6-(6+5+1)]\mod 7=(-6)\mod 7=1.
\end{equation*}
Therefore the Fridays in July 1718 are: \textbf{1, 8, 15, 22} and \textbf{29}.
\end{example}

\begin{example}
Find all Wenesdays in January 1972. \\
\textit{Solution}: \\
$W$ (the code for Wedesday) is 4,\\
$M$ (the code for the January) is 0,\\
$C$ (the code for the Century: 19) is 1, \\
$Y$ (the code for the Year: 72) is 6. \\
Since 1972 is a leap year and the month is January, the date of the first Wednesday is 
\begin{equation*}
D=[4-(0+1+6)+1]\mod 7=(-2)\mod 7=5.
\end{equation*}
Therefore the Wednesdays in January 1972 are: \textbf{5, 12, 19} and \textbf{26}.
\end{example}

\begin{example}
Find all Sundays in February 2080. \\
\textit{Solution}: \\
$W$ (the code for Sunday) is 1,\\
$M$ (the code for the February) is 3,\\
$C$ (the code for the Century: 20) is 0, \\
$Y$ (the code for the Year: 80) is 2. \\
Since 2080 is a leap year and the month is February, the date of the first Sunday is 
\begin{equation*}
D=[1-(3+0+2)+1]\mod 7=(-2)\mod 7=5.
\end{equation*}
Therefore the Sundays in February 2080 are: \textbf{4, 11, 18} and \textbf{25}.
\end{example}

%%%%%%%%%%%%%%%%%%%%%%%%%%%%%%%%%%%%%%%%%%%%%%%%%%%%
\section{Some Concluding Remarks }
%%%%%%%%%%%%%%%%%%%%%%%%%%%%%%%%%%%%%%%%%%%%%%%%%%%%

\begin{enumerate}[(a)]% numbered list
\item The calendar of each year repeats itself entirely after every 28 years. Because of this, in Table \ref{tab:Code-Y}, it suffices to know the codes for the years from 00 to 27. The code for any other year can be obtained from these. Note that the years 00, 28, 56 and 84 have the same codes since they are repeated years. The codes for the other years are obtained as follows. Let $F$, $G$ and $H$ be any year such that $29\leq F \leq 55$, $57\leq G \leq 83$,  and $85\leq H \leq 99$. To get the code for $F$, subtract 28 from it and then pick the code for the corresponding year from 00 to 27.  To get the code for $G$, subtract 56 from it and then pick the code for the corresponding year from 00 to 27. To get the code for $H$, subtract 84 from it and then pick the code for the corresponding year from 00 to 27. 
\item Proposition \ref{prop:Gauss-DF} and Proposition \ref{prop:Nantomah-DF} are applicable only in the Gregorian calendar system. Specifically, they are valid for calculating dates from 15th October 1582 onwards.
\end{enumerate}

%----------------------------------------------------BIBLIOGRAGHY-------------------------------------------------------------
\bibliographystyle{plain}

%------------------------------------------------------------------------------------------------------------------------------------

\end{document}